\documentclass{article}
\usepackage{amssymb,amsmath,theorem,euscript}

\input{epsf}

\newcounter{sec}

\def\sm{\smallskip}

\newcounter{punct}[sec]

\def\punct{\refstepcounter{punct}{\arabic{sec}.\arabic{punct}.  }}

\def\COUNTERS{\addtocounter{sec}{1}
              \setcounter{punct}{0}
          \setcounter{equation}{0}
          \setcounter{theorem}{0}
                  }
\newtheorem{theorem}{Theorem}[sec]

\newtheorem{corollary}[theorem]{Corollary}

\begin{document}

 \def\ov{\overline}
\def\wt{\widetilde}
 \newcommand{\rk}{\mathop {\mathrm {rk}}\nolimits}
\newcommand{\Aut}{\mathop {\mathrm {Aut}}\nolimits}
\newcommand{\Out}{\mathop {\mathrm {Out}}\nolimits}
\renewcommand{\Re}{\mathop {\mathrm {Re}}\nolimits}
\renewcommand{\Im}{\mathop {\mathrm {Im}}\nolimits}
 \newcommand{\tr}{\mathop {\mathrm {tr}}\nolimits}
  \newcommand{\Hom}{\mathop {\mathrm {Hom}}\nolimits}
   \newcommand{\diag}{\mathop {\mathrm {diag}}\nolimits}
   \newcommand{\supp}{\mathop {\mathrm {supp}}\nolimits}
 \newcommand{\im}{\mathop {\mathrm {im}}\nolimits}

\def\Br{\mathrm {Br}}

\def\SL{\mathrm {SL}}
\def\SU{\mathrm {SU}}
\def\GL{\mathrm {GL}}
\def\U{\mathrm U}
\def\OO{\mathrm O}
 \def\Sp{\mathrm {Sp}}
 \def\SO{\mathrm {SO}}
\def\SOS{\mathrm {SO}^*}
 \def\Diff{\mathrm{Diff}(S^1)}
 \def\Vect{\mathfrak{Vect}}
\def\PGL{\mathrm {PGL}}
\def\PU{\mathrm {PU}}
\def\PSL{\mathrm {PSL}}
\def\Symp{\mathrm{Symp}}
\def\End{\mathrm{End}}
\def\Mor{\mathrm{Mor}}
\def\Aut{\mathrm{Aut}}
 \def\PB{\mathrm{PB}}
 \def\cA{\mathcal A}
\def\cB{\mathcal B}
\def\cC{\mathcal C}
\def\cD{\mathcal D}
\def\cE{\mathcal E}
\def\cF{\mathcal F}
\def\cG{\mathcal G}
\def\cH{\mathcal H}
\def\cJ{\mathcal J}
\def\cI{\mathcal I}
\def\cK{\mathcal K}
 \def\cL{\mathcal L}
\def\cM{\mathcal M}
\def\cN{\mathcal N}
 \def\cO{\mathcal O}
\def\cP{\mathcal P}
\def\cQ{\mathcal Q}
\def\cR{\mathcal R}
\def\cS{\mathcal S}
\def\cT{\mathcal T}
\def\cU{\mathcal U}
\def\cV{\mathcal V}
 \def\cW{\mathcal W}
\def\cX{\mathcal X}
 \def\cY{\mathcal Y}
 \def\cZ{\mathcal Z}
\def\0{{\ov 0}}
 \def\1{{\ov 1}}
 \def\frA{\mathfrak A}
 \def\frB{\mathfrak B}
\def\frC{\mathfrak C}
\def\frD{\mathfrak D}
\def\frE{\mathfrak E}
\def\frF{\mathfrak F}
\def\frG{\mathfrak G}
\def\frH{\mathfrak H}
\def\frI{\mathfrak I}
 \def\frJ{\mathfrak J}
 \def\frK{\mathfrak K}
 \def\frL{\mathfrak L}
\def\frM{\mathfrak M}
 \def\frN{\mathfrak N} \def\frO{\mathfrak O} \def\frP{\mathfrak P} \def\frQ{\mathfrak Q} \def\frR{\mathfrak R}
 \def\frS{\mathfrak S} \def\frT{\mathfrak T} \def\frU{\mathfrak U} \def\frV{\mathfrak V} \def\frW{\mathfrak W}
 \def\frX{\mathfrak X} \def\frY{\mathfrak Y} \def\frZ{\mathfrak Z} \def\fra{\mathfrak a} \def\frb{\mathfrak b}
 \def\frc{\mathfrak c} \def\frd{\mathfrak d} \def\fre{\mathfrak e} \def\frf{\mathfrak f} \def\frg{\mathfrak g}
 \def\frh{\mathfrak h} \def\fri{\mathfrak i} \def\frj{\mathfrak j} \def\frk{\mathfrak k} \def\frl{\mathfrak l}
 \def\frm{\mathfrak m} \def\frn{\mathfrak n} \def\fro{\mathfrak o} \def\frp{\mathfrak p} \def\frq{\mathfrak q}
 \def\frr{\mathfrak r} \def\frs{\mathfrak s} \def\frt{\mathfrak t} \def\fru{\mathfrak u} \def\frv{\mathfrak v}
 \def\frw{\mathfrak w} \def\frx{\mathfrak x} \def\fry{\mathfrak y} \def\frz{\mathfrak z} \def\frsp{\mathfrak{sp}}
 \def\bfa{\mathbf a} \def\bfb{\mathbf b} \def\bfc{\mathbf c} \def\bfd{\mathbf d} \def\bfe{\mathbf e} \def\bff{\mathbf f}
 \def\bfg{\mathbf g} \def\bfh{\mathbf h} \def\bfi{\mathbf i} \def\bfj{\mathbf j} \def\bfk{\mathbf k} \def\bfl{\mathbf l}
 \def\bfm{\mathbf m} \def\bfn{\mathbf n} \def\bfo{\mathbf o} \def\bfp{\mathbf p} \def\bfq{\mathbf q} \def\bfr{\mathbf r}
 \def\bfs{\mathbf s} \def\bft{\mathbf t} \def\bfu{\mathbf u} \def\bfv{\mathbf v} \def\bfw{\mathbf w} \def\bfx{\mathbf x}
 \def\bfy{\mathbf y} \def\bfz{\mathbf z} \def\bfA{\mathbf A} \def\bfB{\mathbf B} \def\bfC{\mathbf C} \def\bfD{\mathbf D}
 \def\bfE{\mathbf E} \def\bfF{\mathbf F} \def\bfG{\mathbf G} \def\bfH{\mathbf H} \def\bfI{\mathbf I} \def\bfJ{\mathbf J}
 \def\bfK{\mathbf K} \def\bfL{\mathbf L} \def\bfM{\mathbf M} \def\bfN{\mathbf N} \def\bfO{\mathbf O} \def\bfP{\mathbf P}
 \def\bfQ{\mathbf Q} \def\bfR{\mathbf R} \def\bfS{\mathbf S} \def\bfT{\mathbf T} \def\bfU{\mathbf U} \def\bfV{\mathbf V}
 \def\bfW{\mathbf W} \def\bfX{\mathbf X} \def\bfY{\mathbf Y} \def\bfZ{\mathbf Z} \def\bfw{\mathbf w}
 \def\R {{\mathbb R }} \def\C {{\mathbb C }} \def\Z{{\mathbb Z}} \def\H{{\mathbb H}} \def\K{{\mathbb K}}
 \def\N{{\mathbb N}} \def\Q{{\mathbb Q}} \def\A{{\mathbb A}} \def\T{\mathbb T} \def\P{\mathbb P} \def\G{\mathbb G}
 \def\bbA{\mathbb A} \def\bbB{\mathbb B} \def\bbD{\mathbb D} \def\bbE{\mathbb E} \def\bbF{\mathbb F} \def\bbG{\mathbb G}
 \def\bbI{\mathbb I} \def\bbJ{\mathbb J} \def\bbL{\mathbb L} \def\bbM{\mathbb M} \def\bbN{\mathbb N} \def\bbO{\mathbb O}
 \def\bbP{\mathbb P} \def\bbQ{\mathbb Q} \def\bbS{\mathbb S} \def\bbT{\mathbb T} \def\bbU{\mathbb U} \def\bbV{\mathbb V}
 \def\bbW{\mathbb W} \def\bbX{\mathbb X} \def\bbY{\mathbb Y} \def\kappa{\varkappa} \def\epsilon{\varepsilon}
 \def\phi{\varphi} \def\le{\leqslant} \def\ge{\geqslant}

\def\UU{\bbU}
\def\Mat{\mathrm{Mat}}
\def\tto{\rightrightarrows}

\def\Gms{\mathrm {Gms}}
\def\Ams{\mathrm {Ams}}
\def\Isom{\mathrm {Isom}}

\def\Gr{\mathrm{Gr}}

\def\graph{\mathrm{graph}}

\def\O{\mathrm{O}}

\def\la{\langle}
\def\ra{\rangle}


 \def\ov{\overline}
\def\wt{\widetilde}

\renewcommand{\Re}{\mathop {\mathrm {Re}}\nolimits}
\def\Br{\mathrm {Br}}

  \def\ASp{\mathrm {ASp}}
 \def\Isom{\mathrm {Isom}}
\def\SL{\mathrm {SL}}
\def\SU{\mathrm {SU}}
\def\GL{\mathrm {GL}}
\def\U{\mathrm U}
\def\OO{\mathrm O}
 \def\Sp{\mathrm {Sp}}
 \def\SO{\mathrm {SO}}
\def\SOS{\mathrm {SO}^*}
 \def\Vect{\mathfrak{Vect}}
\def\PGL{\mathrm {PGL}}
\def\PU{\mathrm {PU}}
\def\PSL{\mathrm {PSL}}
\def\Symp{\mathrm{Symp}}
\def\End{\mathrm{End}}
\def\Mor{\mathrm{Mor}}
\def\Aut{\mathrm{Aut}}
 \def\PB{\mathrm{PB}}
 \def\cA{\mathcal A}
\def\cB{\mathcal B}
\def\cC{\mathcal C}
\def\cD{\mathcal D}
\def\cE{\mathcal E}
\def\cF{\mathcal F}
\def\cG{\mathcal G}
\def\cH{\mathcal H}
\def\cJ{\mathcal J}
\def\cI{\mathcal I}
\def\cK{\mathcal K}
 \def\cL{\mathcal L}
\def\cM{\mathcal M}
\def\cN{\mathcal N}
 \def\cO{\mathcal O}
\def\cP{\mathcal P}
\def\cQ{\mathcal Q}
\def\cR{\mathcal R}
\def\cS{\mathcal S}
\def\cT{\mathcal T}
\def\cU{\mathcal U}
\def\cV{\mathcal V}
 \def\cW{\mathcal W}
\def\cX{\mathcal X}
 \def\cY{\mathcal Y}
 \def\cZ{\mathcal Z}
\def\0{{\ov 0}}
 \def\1{{\ov 1}}
 \def\frA{\mathfrak A}
 \def\frB{\mathfrak B}
\def\frC{\mathfrak C}
\def\frD{\mathfrak D}
\def\frE{\mathfrak E}
\def\frF{\mathfrak F}
\def\frG{\mathfrak G}
\def\frH{\mathfrak H}
\def\frI{\mathfrak I}
 \def\frJ{\mathfrak J}
 \def\frK{\mathfrak K}
 \def\frL{\mathfrak L}
\def\frM{\mathfrak M}
 \def\frN{\mathfrak N} \def\frO{\mathfrak O} \def\frP{\mathfrak P} \def\frQ{\mathfrak Q} \def\frR{\mathfrak R}
 \def\frS{\mathfrak S} \def\frT{\mathfrak T} \def\frU{\mathfrak U} \def\frV{\mathfrak V} \def\frW{\mathfrak W}
 \def\frX{\mathfrak X} \def\frY{\mathfrak Y} \def\frZ{\mathfrak Z} \def\fra{\mathfrak a} \def\frb{\mathfrak b}
 \def\frc{\mathfrak c} \def\frd{\mathfrak d} \def\fre{\mathfrak e} \def\frf{\mathfrak f} \def\frg{\mathfrak g}
 \def\frh{\mathfrak h} \def\fri{\mathfrak i} \def\frj{\mathfrak j} \def\frk{\mathfrak k} \def\frl{\mathfrak l}
 \def\frm{\mathfrak m} \def\frn{\mathfrak n} \def\fro{\mathfrak o} \def\frp{\mathfrak p} \def\frq{\mathfrak q}
 \def\frr{\mathfrak r} \def\frs{\mathfrak s} \def\frt{\mathfrak t} \def\fru{\mathfrak u} \def\frv{\mathfrak v}
 \def\frw{\mathfrak w} \def\frx{\mathfrak x} \def\fry{\mathfrak y} \def\frz{\mathfrak z} \def\frsp{\mathfrak{sp}}
 \def\bfa{\mathbf a} \def\bfb{\mathbf b} \def\bfc{\mathbf c} \def\bfd{\mathbf d} \def\bfe{\mathbf e} \def\bff{\mathbf f}
 \def\bfg{\mathbf g} \def\bfh{\mathbf h} \def\bfi{\mathbf i} \def\bfj{\mathbf j} \def\bfk{\mathbf k} \def\bfl{\mathbf l}
 \def\bfm{\mathbf m} \def\bfn{\mathbf n} \def\bfo{\mathbf o} \def\bfp{\mathbf p} \def\bfq{\mathbf q} \def\bfr{\mathbf r}
 \def\bfs{\mathbf s} \def\bft{\mathbf t} \def\bfu{\mathbf u} \def\bfv{\mathbf v} \def\bfw{\mathbf w} \def\bfx{\mathbf x}
 \def\bfy{\mathbf y} \def\bfz{\mathbf z} \def\bfA{\mathbf A} \def\bfB{\mathbf B} \def\bfC{\mathbf C} \def\bfD{\mathbf D}
 \def\bfE{\mathbf E} \def\bfF{\mathbf F} \def\bfG{\mathbf G} \def\bfH{\mathbf H} \def\bfI{\mathbf I} \def\bfJ{\mathbf J}
 \def\bfK{\mathbf K} \def\bfL{\mathbf L} \def\bfM{\mathbf M} \def\bfN{\mathbf N} \def\bfO{\mathbf O} \def\bfP{\mathbf P}
 \def\bfQ{\mathbf Q} \def\bfR{\mathbf R} \def\bfS{\mathbf S} \def\bfT{\mathbf T} \def\bfU{\mathbf U} \def\bfV{\mathbf V}
 \def\bfW{\mathbf W} \def\bfX{\mathbf X} \def\bfY{\mathbf Y} \def\bfZ{\mathbf Z} \def\bfw{\mathbf w}
 \def\R {{\mathbb R }} \def\C {{\mathbb C }} \def\Z{{\mathbb Z}} \def\H{{\mathbb H}} \def\K{{\mathbb K}}
 \def\N{{\mathbb N}} \def\Q{{\mathbb Q}} \def\A{{\mathbb A}} \def\T{\mathbb T} \def\P{\mathbb P} \def\G{\mathbb G}
 \def\bbA{\mathbb A} \def\bbB{\mathbb B} \def\bbD{\mathbb D} \def\bbE{\mathbb E} \def\bbF{\mathbb F} \def\bbG{\mathbb G}
 \def\bbI{\mathbb I} \def\bbJ{\mathbb J} \def\bbL{\mathbb L} \def\bbM{\mathbb M} \def\bbN{\mathbb N} \def\bbO{\mathbb O}
 \def\bbP{\mathbb P} \def\bbQ{\mathbb Q} \def\bbS{\mathbb S} \def\bbT{\mathbb T} \def\bbU{\mathbb U} \def\bbV{\mathbb V}
 \def\bbW{\mathbb W} \def\bbX{\mathbb X} \def\bbY{\mathbb Y} \def\kappa{\varkappa} \def\epsilon{\varepsilon}
 \def\phi{\varphi} \def\le{\leqslant} \def\ge{\geqslant}

\def\UU{\bbU}
\def\Mat{\mathrm{Mat}}
\def\tto{\rightrightarrows}

\def\Gr{\mathrm{Gr}}

\def\graph{\mathrm{graph}}

\def\Vir{\mathfrak{Vir}}

\def\O{\mathrm{O}}

\def\la{\langle}
\def\ra{\rangle}

\def\wre{\divideontimes}

\begin{center}
\bf\large
The group of diffeomorphisms of the circle:
 reproducing kernels and analogs of spherical functions 
 
 \bigskip
 
 \sc Yury Neretin%
 \footnote{Supported by the grants FWF, projects P25142, P28421.}
\end{center}

{\small The group $\Diff$ of diffeomorphisms of the circle is an infinite dimensional analog of the
real semisimple Lie groups $\U(p,q)$, $\Sp(2n,\R)$, $\SO^*(2n)$; the space $\Xi$ of univalent functions
is an analog of the corresponding
 classical complex Cartan domains. We present explicit formulas for realizations
of highest weight representations of $\Diff$ in the space of holomorphic functionals on $\Xi$, reproducing kernels
on $\Xi$ determining inner products, and expressions ('canonical cocycles') replacing spherical functions.}

\section{Introduction}

\COUNTERS

{\bf\punct The purpose of the paper.}
The group $\Diff$ of orientation preserving diffeomorphisms of the circle has
unitary
projective highest weight  representations.
There is a  well-developed representation theory
of unitary highest weight representations of real semi-simple Lie groups (see, e.g., an elementary introduction in \cite{Ner-gauss}, Chapter 7, and  further references 
in this book).
In particular, this theory includes realizations of highest weight representations in spaces of holomorphic (vector-valued) functions
on Hermitian symmetric spaces, reproducing kernels, Berezin--Wallach sets, Olshanski semigroups,
Berezin--Guichardet--Wigner formulas for central extensions.
All  these phenomena exist for the group $\Diff$.

The purpose of this paper is to present details of a strange calculation sketched  at the end of my paper
\cite{Ner-sbornik}, 1989.
It  also was included to  my  thesis \cite{Ner-disser}, 1992.  Later the calculation and its result never
were exposed or repeated.

Let us consider a unitary highest weight representation $\rho$ of  the group $\Diff$. Generally, it is 
a projective representation. Operators $\rho(g)$ are determined up to  constant factors.
 We 
 normalize them by the condition
$$
\la \rho(g) v, v\ra=1,
$$
where $v$ is a highest weight vector%
\footnote{Since the representation $\rho$ is projective, the operators $\rho(g)$ are defined up to scalar factors. If $\la \rho(g) v, v\ra\ne 0$, then we can write corrected operators
	\newline
	$\phantom{.}\qquad\qquad\wt\rho(g):=\bigl(\la \rho(g) v, v\ra\bigr)^{-1} \rho(g)$.
	\newline
The property  $\la \rho(g) v, v\ra\ne 0$ holds for all $g\in\Diff$.
	This follows from explicit calculations given below. Also, it is possible a priory proof from
	formula (\ref{eq:power}).
 Of course, the new operators
	$\wt \rho(g)$ are not unitary.}.
  Then we have
$$
\rho(g_1)\rho(g_2)=c(g_1,g_2)\, \rho(g_1 g_2).
$$
where the expression $c(g_1,g_2)$ (the {\it canonical cocycle}) is canonically determined
by the representation $\rho$. 
In this paper, we derive a formula for $c(g_1,g_2)$.

\sm

{\sc Remark.} Let $G$ be a semisimple group, $K$  a maximal compact subgroup
and $v$  a $K$-fixed vector, let $\phi(g)=\la\rho(g) v, v \ra$ be the spherical function.
Then
$$
\qquad\qquad\qquad\qquad\qquad
c(g_1,g_2)=\frac{\phi(g_1g_2)}{\phi(g_1)\,\phi(g_2)}
.\qquad\qquad\qquad\qquad\qquad\quad
$$
In our case, spherical functions are not well-defined but the 'canonical cocycles',
which are hybrids of central extensions and spherical functions, exist. \hfill $\boxtimes$

\sm

A formula for canonical cocycles implies some  automatic corollaries:
 we can write explicit formulas for realizations
of highest weight representations of $\Diff$ on the space univalent functions
and for invariant reproducing kernels on the space of holomorphic functionals on
the space of univalent functions.

\sm

{\bf\punct Virasoro algebra and highest weight modules.%
\label{ss:virasoro}} For details, see, e.g., \cite{Ner-affine}.
Denote by $\Vect_\R$ and $\Vect_\C$ respectively the Lie algebras of real (respectively
complex) vector fields on the circle. Choosing a basis
$L_n:=e^{in\phi}\partial/i\partial\phi$ in $\Vect_\C$, we get the commutation relations
$$
 [L_n,L_m]=(m-n) L_{m+n}.
 $$

 Recall that the {\it Virasoro algebra} $\Vir$ is a Lie algebra with a basis $L_n$,  $\zeta$, where $n$ ranges in $\Z$, 
 and commutation relations
 $$
 [L_n,L_m]=(m-n) L_{m+n}+\frac{n^3-n}{12} \delta_{n+m,0} \zeta,\qquad [L_n,\zeta]=0.
 $$

Let $h$, $c\in \C$. A {\it module with the   highest weight} $(h,c)$  over $\Vir$ is a module containing a vector $v$
 such that 
 
 \sm
 
 1) $L_0 v=h v$, $\zeta v= cv$;
 
 \sm
 
 2) $L_{-n} v=$ for $n<0$;
 
 \sm
 
 3) the vector $v$ is cyclic.
 
 \sm
 
 There exists a unique irreducible  highest weight module   $L(h,c)$ with a given highest weight $(h,c)$, also
 there is a universal highest weight module $M(h,c)$ (a {\it Verma module}), such that any module with highest weight $(h,c)$
 is a quotient of the Verma module $M(h,c)$. In particular, if $M(h,c)$
 is irreducible, we have $M(h,c)=L(h,c)$.
 
 By  \cite{Kac}, \cite{FF1}, a Verma module $M(h,c)$ is reducible  if and only if  $(h,c)\in \C^2$
 satisfies at least one equation of the form
 \begin{multline*}
 \Bigl(h-\frac 1{24}(\alpha^2-1)(c-13)+\frac 12 (\alpha^2-1)  \Bigr)
 \Bigl(h-\frac 1{24}(\beta^2-1)(c-13)+\frac 12 (\beta^2-1)  \Bigr)
 +\\+
 \frac 1{16} (\alpha^2-\beta^2)^2=0,\qquad \text{where $\alpha$, $\beta$ range in $\N$.}
 \end{multline*}

  \sm
  
 If $M(h,c)$ is reducible, its composition series is finite or countable, it is described in \cite{FF2},  \cite{FF3}.

\sm

{\bf\punct Unitarizability.%
\label{ss:initarizability}} A module $L(h,c)$ is called {\it unitarizable} if it admits a positive definite inner
product such that $L_n=-L_{-n}^*$.  A module $L(h,c)$ is unitarizable iff  $(h,c)\in \R^2$ satisfies
one of the following conditions:
\begin{align}
&\text{1. (continulous series)}\,\,
 h\ge 0,\qquad c\ge 1;
\label{eq:unitarity-1}
\\
\!\!\!\!\!\!\!\!\!\!
&\text{2. (discrete series)}\,\,
c=1-\frac 6{p(p+1)},\qquad h=\frac{\bigl(\alpha p-\beta(p+1)\bigr)^2-1}{4p(p+1)},
\label{eq:unitarity-2}
\end{align}
where $\alpha$, $\beta$, $p\in \Z$, $p\ge 2$, $1\le \alpha \le p$, $1\le \beta \le p-1$.

 \sm
 
 For sufficiency of these conditions, see \cite{GKO}, for necessity, see \cite{FQS}, \cite{Ner-disser83}, \cite{Ner-GB}.

 \sm
 
 {\sc Remark.} 1) For $h>0$, $c>1$ we have $L(h,c)=M(h,c)$.
 For $(h,c)$ of the form (\ref{eq:unitarity-2}) modules $M(h,c)$ have countable
 composition series.
 
 \sm
 
2) There are simple explicit constructions for representations 
 $L(h,c)$ for $c\ge 1$, $h\ge(c-1)/24$. The present paper is based on a such construction,
 for another (fermionic) construction, see \cite{Ner-disser83}, \cite{Ner-book}, Sect. 7.3).
 
 \sm
 
3) There is an explicit construction for representations (\ref{eq:unitarity-2})
 with $c=1/2$ (possible values of $h$ are 0, $1/16$, $1/2$), see
 \cite{Ner-1983}, \cite{Ner-book}, Sect. 7.3. The module $L(0,0)$ is the
 trivial one-dimensional module.
 
 \sm
 
4)   As far as I know, transparent constructions (that visualize
the action of the algebra/group, the space of
  representation and the inner product) for the other points of discrete series
  and for the domain $0\le h< (c-1)/24$ are unknown.
  \hfill $\boxtimes$
  
  \sm
 
 {\bf\punct Standard boson realizations.%
 \label{ss:boson-realization}}
 Consider the space $\cF$ of polynomials of variables $z_1$, $z_2$, \dots.
 Define the {\it creation and annihilation operators} $a_n$, where $n$ ranges
 in $\Z\setminus 0$, by
 $$
 a_nf(z)=
 \begin{cases}\sqrt n z_n f(z),\qquad\text{for $n>0$};
 \\
  \sqrt {(-n)}\cdot\frac\partial {\partial z_n} f(z),
 \qquad \text{for $n<0$}
 .
 \end{cases}
 $$
 
 Next, we write formulas for representations of the Virasoro
 algebra in $\cF$. Fix the parameters $\alpha$, $\beta\in \C$.
 For $n\ne0$ we set 
 \begin{eqnarray*}
 L_n&:=&\frac 12\sum_{k,l:\,k+l=n} a_k a_l+ (\alpha+in\beta)a_n;
 \\
 L_0&:=& \sum_{n>0} a_n a_{-n}+\frac 12 (\alpha^2+\beta^2);
 \\
 \zeta&:=& 1+12\beta^2.
 \end{eqnarray*}
 
Then we have $L_{-n} 1=0$ for $n<0$ and
$$
L_0\cdot 1=\frac1 2(\alpha^2+\beta^2),\qquad  \zeta\cdot 1:= 1+12\beta^2.
$$
In this way, we get a representation of $\Vir$, whose Jordan--H\"older series coincides with the
Jordan--H\"older series of $M(h,c)$ with 
$$
h=\frac 12(\alpha^2+\beta^2), \qquad c=1+12\beta^2.
$$
For points $(\alpha,\beta)$ in a general position we get a Verma module.

\sm

{\sc Remarks.} 1) For formulas for singular vectors, see \cite{Mim}.

\sm

2) For $\alpha$ and $\beta$ in a general position representations 
with parameters $(\alpha,\beta)$ and $(\alpha,-\beta)$ are equivalent.
As far as I know explicit formula for the intertwining operator remains to
be uknown. \hfill $\boxtimes$

\sm

{\bf \punct The group of diffeomorphisms of the circle.%
\label{ss:Diff}}
Denote by $\Diff$ the group of smooth orientation preserving diffeomorphisms of the circle.
Its Lie algebra is $\Vect_\R$.

Any unitarizable highest weight representation of $\Vir$  can be integrated to a projective unitary representation
of $\Diff$, see \cite{GW}. A module $L(h,c)$ with an arbitrary highest weight
$(h,c)\in\C^2$ can be integrated to a
projective representation
of $\Diff$ by bounded operators in a Frechet space, see \cite{Ner-spinor}.
This follows from the universal fermionic construction, see \cite{Ner-book}, Sect. VII.3. 

\sm

{\bf \punct The welding.%
\label{ss:welding}} Denote by $\ov \C=\C\cup \infty$ the Riemann sphere.
Denote by $D_+$ the disk $|z|\le 1$ on $\ov\C$, by $D_-$ the disk $|z|\ge 1$ on $\ov\C$.
Denote by $D_\pm^\circ$ their interiors. Denote by $S^1$ the circle $|z|=1$.

We say that a function $f:D_\pm \to\ov \C$
is {\it univalent up to the boundary} if $f$ is an embedding, 
which is  holomorphic  in the interior of the disk and is smooth up to the boundary.

Let $\gamma\in\Diff$. Let us glue the disks $D_+$ and $D_-$ identifying the points
$z\in S^1\subset D_+$ with $\gamma(z)\in S^1\subset D_-$. We get a two-dimensional real manifold with
complex structure on the images of $D_+$ and $D_-$. According \cite{AB} this structure
admits a unique extension to the separating contour (in fact, the smoothness of $\gamma$ is redundant, for a minimal
condition, see \cite{AB}). Thus we get a one-dimensional complex manifold, i.e., a Riemann sphere $\ov\C$,
and a pair of univalent maps $D_\pm\to\ov \C$.

Conversely, consider a Riemann sphere $\ov\C$ and  a pair of maps $p_-:D_-\to \ov\C$,
$p_+:D_+\to \ov\C$
univalent up to a boundary such that
$$
p_-(D_-^\circ)\cap p_+(D_+^\circ)=\varnothing \qquad p_-(D_-)\cup p_+(D_+)=\ov\C
.$$
Then we have a diffeomorphism $p_-^{-1}\circ p_+:S^1\to S^1$.

\sm

{\bf \punct The semigroup $\Gamma$.%
\label{ss:Gamma}}
An element of the semigroup%
\footnote{This semigroup is an analog of Olshanski semigroups,
	see, e.g., \cite{Ner-gauss}, Add. to Chapter 3.}
$\Gamma$ is a triple $\frR:=(\cR, r_+,r_-)$, where
$\cR$ is the Riemann surface (a one-dimensional complex manifold)
equivalent to the Riemann sphere $\ov\C$, 
$$
r_-:D_-\to \cR, \qquad r_+:D_+\to \cR
$$
are
univalent up to the boundary, and 
$$
r_-(D_-)\cap r_+(D_+)=\varnothing.
$$
Two triples  $(\cR, r_+,r_-)$ and $(\cR', r_+',r_-')$ are equivalent
if there is a biholomorphic map $\phi:\cR\to \cR'$ such that
$r'_\pm=\phi \circ r_\pm$.

Define a product $\frR=\frP\frQ$ of two elements
$\frP:=(\cP, p_+,p_-)$, $\frQ:=(\cQ, q_+,q_-)$.
We take surfaces (closed disks) 
$$\cP\setminus p_-(D_-^\circ),\qquad
\cQ\setminus q_+(D_+^\circ)$$
 and glue them together 
identifying points $p_-(z)\in \cP\setminus p_-(D_-^\circ)$
and $q_+(z)\in \cQ\setminus q_+(D_+^\circ)$, where $z$ ranges in $S^1$.
We get a Riemann sphere and two univalent functions, i.e., an element of $\Gamma$.

It is clear that diffeomorphisms $\gamma\in\Diff$ can be regarded as 
limit points of $\Gamma$.
We also define the semigroup 
$$\ov \Gamma=\Gamma\cup\Diff.$$

This semigroup was discovered in \cite{Ner-semi}, \cite{Seg}; for details, see \cite{Ner-sbornik}, \cite{Ner-book}, Sect VII.4.

\sm

There are the following statements about representations of  $\Gamma$:

\sm

a) Any module $L(h,c)$ can be integrated to a projective holomorphic representation 
of the semigroup $\ov\Gamma$ by bounded operators in a Frechet space.

\sm

b)  Any unitarizable module
$L(h,c)$ can be integrated to a representation of $\ov \Gamma$, which is 
holomorphic on $\Gamma$ and unitary up to scalars on $\Diff$, see \cite{Ner-sbornik}, \cite{Ner-book}, Sect. 7.4.

\sm

c) Standard boson represantations of Virasoro algebra (see Subsect. \ref{ss:boson-realization})
can be integrated  to holomorphic representations of $\Gamma$ by bounded operators,
see \cite{Ner-sbornik}. Our calculation is based on explicit formulas for such representations.

\sm

{\bf \punct  The complex domain $\Diff/\bbT$.%
\label{ss:Xi}}
Next, consider the space $\Xi$, whose points are pairs
$(\cS,s,\sigma)$, where $\cS$ is a Riemann surface equivalent to Riemannian
sphere and $s:D_+\to \cS$ is a map univalent up to the boundary,
and $\sigma$ is a point in $\cS\setminus s(D_+)$.
The semigroup $\ov\Gamma$ acts on $\Xi$
in the following way.
Let $\frP:=(\cP, p_+,p_-)\in\Gamma$, $\frS=(\cS,s)$.
We glue $\cP\setminus p_-(D_-^0)$ and $\cS\setminus s(D_+^\circ)$ identifying points
$p_-(z)\in \cP\setminus p_-(D_-^0)$ and $s(z)\in \cS\setminus s(D_+^\circ)$,
where $z$ ranges $S^1$. We get a Riemann sphere, an univalent map from $D_+$
to the Riemann sphere, and a distinguished point.

Below we assume
$$
\cS=\ov \C,\quad s(0)=0,\qquad \sigma=\infty.
$$

{\sc Remark.} The space $\wt\Xi$ of univalent functions $s(z)=z+a_2z^2+\dots$
from $D_+^\circ$ to $\C$ was a subject of a wide and interesting literature
(see books \cite{Gol}, \cite{Duren}). 
In 1907 Koebe showed that this space is compact, in next 80 years extrema of 
numerous functionals on $\wt\Xi$ were explicitly evaluated. Infortunately,
this science declared the Biebarbach conjecture as the main purpose.
The proof of the conjecture (De Branges, 1985)  implied a decay of interest 
to the subject.
Unforunately, this happened  before $\Xi$ and $\wt\Xi$ became a topic of  the representation theory
 \cite{Kir1}, \cite{Ner-semi}, \cite{BR}, \cite{Ner-sbornik} in 1986-1989.
\hfill $\boxtimes$

\sm

{\bf\punct The cocycles $\lambda$ and $\mu$.%
\label{ss:lambda-mu}}
Let 
$$
r_-:D_-\to \ov\C\setminus 0, \qquad p_+:D_+\to \ov\C\setminus\infty
$$
 be functions univalent up to the boundary.
We assume that $r_-(\infty)=\infty$, $p_+(0)=0$ (there are no extra conditions).

Let us take a function $r^+:D_+\to \ov \C$ such that $r^+(0)=0$ and $(\ov\C, r_-, r^+)$
determines an element $\Diff$, denote it%
\footnote{There is a one-parametric family of such diffeomorphisms} by $\gamma_r$. Take a function
$p^-: D_-\to\ov\C$ such that $(\ov\C, p_+,p^-)$ determines an element of $\Diff$, denote it
by $\gamma_p$.
Consider the product $\gamma_r \gamma_p$ and take the corresponding triple
$$\gamma_r\gamma_p=(\cQ,q_+,q_-) .$$ 
We can  assume $\cQ=\ov\C$, $q_+(0)=0$, $q_-(\infty)=\infty$.

We define two functions%
\footnote{The formulas in \cite{Ner-sbornik} contain
	a mistake in rational coefficients in the front of the integrals. This
	happened in 
 maniplulations with 'logarithmic forms' $f(z)+\ln dz$.
	 We have $\ln dz=i\phi+\ln d\phi$,  the summand $i\phi$ \cite{Ner-sbornik} was lost.} 
\begin{equation}
\lambda(r_-,p_+)=\int_{|z|=1} \Bigl( \ln \frac{r_-(z)}{z}\,\, d\ln \frac{q_+(\gamma_r(z))}{r^+(\gamma_r(z))}-
\ln \frac{p_+(\gamma_p(z))}{\gamma_p(z)} \,\, d\ln \frac{q_-(z)}{p^-(z)} \Bigr).
\label{eq:lambda}
\end{equation}
\begin{multline}
\mu (r_-,p_+)=\\=
\frac 1{24}\int_{|z|=1} \Bigl( -\ln r_-(z)'\,\, d\ln \frac{q_+'(\gamma_r(z))}{(r^+)'(\gamma_r(z))}+
\ln p_+'(\gamma_p(z))\, \,d\ln\frac{q_-'(z)}{(p^-)'(z)}\Bigr)
.
\label{eq:mu}
\end{multline}


{\bf \punct Canonical cocycles.%
\label{ss:cocycles}}
Consider a highest weight representation $L(h,c)$ of $\Vir$. Denote
by $\cN^{h,c}$ the corresponding representation
of the semigroup $\Gamma$. Since our representation is projective, operators
$\cN^{h,c}(\frR)$ are defined up to  scalar factors.

Denote by $v$ a highest weight vector in the completed  $L(h,c)$,
it is determined up to a scalar factor. The projection operator $\pi$ to the line $\C v$
is determined canonically (its  kernel is a sum of all weight subspaces with weight different from $h$).
The condition 
$$
\pi \left(\cN^{h,c}(\frR) v\right)= v
$$
determines an operator $\cN^{h,c}(\frR)$ canonically.
Now we have
\begin{equation}
\cN^{h,c}(\frR)\, \cN^{h,c}(\frP)= \kappa^{h,c} (\frR, \frP)\, \cN^{h,c}(\frR\circ \frP)
,
\label{eq:kappa-hc-1}
\end{equation}
where a constant $\kappa^{h,c} (\frR, \frP)\in \C$ is a canonically defined.

\begin{theorem}
 \begin{equation}
  \kappa^{h,c}(\frR, \frP)= \exp \bigl\{h \lambda(r_-, p_+)\bigr\} \cdot\exp\{c \mu(r_-, p_+)\}.
  \label{eq:kappa-hc-2}
 \end{equation}
\end{theorem}

{\sc Remarks.} 1) By construction, the functions $\lambda$, 
$\mu:\Diff\times\Diff\to \C$ determine $\C$-central extensions
of $\Diff$. In other words they represent classes of the secon group cohomologies 
$H^2(\Diff,\C)$.

\sm

2) The group of {\it continuous} cohomologies $H^2(\Diff,\R)$ is generated by two
cocyles (see \cite{Fuks}, Theorem 3.4.4). The first one is the  {\it Bott cocycle}
$$
c_1(\gamma_1,\gamma_2)=
\frac 1 2
\int_{|z|=1} \ln \bigl(\gamma'_1(\gamma_2(z))\bigr)\,d\ln \gamma_2'(z).
$$
To define the second cocycle, we write the multivalued expression
$$
c_2(\gamma_1,\gamma_2)=
\ln(\gamma_1(\gamma_2(z)))-\ln(\gamma_2(z))-\ln \gamma_1(z)
$$
and choose its continuous branch on $\Diff\times \Diff$ such that $c_1(e,e)=0$
(this cocycle can be reduced to a $\Z$-cocycle determining the universal covering
group 
of  $\Diff$). The cocycles $\mu$, $\lambda$  are equivalent to
$c_1$, $c_2$ in the group $H^2(\Diff,\C)$.

\sm

 3) I never met a theorem that all measurable $\R/2\pi\Z$-central extensions of
$\Diff$ are reduced to cocycles $c_1$, $c_2$. In our context, cocycles are continuous 
by construction exposed below, see (\ref{eq:1-1-1}).  However, the question has some interest
from the point of view of representation theory. 
\hfill $\boxtimes$

\sm

{\bf\punct  Realization of highest weight representations in
	the space of holomorphic functionals  on the domain $\Xi$.%
\label{ss:action-Xi}}

\begin{corollary}
	\label{cor:1}
 The formula 
 \begin{equation}
  \rho_{h,c} (\frQ)\, f(\frS)=f (\frS \frQ) \exp \bigl\{h\, \lambda(q_-, s)\bigr\}
   \cdot\exp\{c\, \mu(q_-,s)\}
  \label{eq:sections}
 \end{equation}
determines a representation of $\ov\Gamma$ in the space of holomorphic functionals
on $\Xi$. For $(h,c)$ in a general position {\rm(}if  the Verma module  $M(h,c)$ is irreducible {\rm)}, the
corresponding
module  over the Virasoro algebra is $M(h,c)$.
\end{corollary}


{\sc Remark.} Formulas for the underlying action of the Lie algebra $\Vect$ on $\Xi$ were firstly obtained in
by Scheffer \cite{Sch}, 1948, and rediscovered in \cite{KYu}.
The action  of $\Vir$ corresponding (\ref{eq:sections}) was present  in \cite{Ner-semi}
without a proof, for a proof, see \cite{Kir2}. Considering an action in the space of polynomials
in Taylor coefficients of $p(z)$, we get a module dual to the Verma module.
Some modifications of formulas
 are contained in \cite{AN}.
 \hfill $\boxtimes$
 
 \sm
 
{\sc Remark.} Of course, we can invent many different definitions of holomorphic functions on $\Xi$
 and of spaces of holomorphic functions. Our statement is almost independent on this.
 To avoid a scholastic discussion, we note that  a construction \cite{Ner-sbornik}, Subs. 4.12
 gives an operator sending   Fock space
 to
 the space of functions on $\Xi$.
  \hfill $\boxtimes$

 \sm
 
 {\bf \punct Reproducing kernel.%
 \label{ss:reproducing-kernel}}
On Hilbert spaces determined by reproducing kernels see, e.g., \cite{Ner-gauss}, Sect. 7.1.
 
 \begin{corollary}
 	\label{cor:2}
  Let $L(h,c)$ be unitarizable. 
  Then the 
  kernel
  \begin{equation}
K^{h,c}(r,p):=\exp \bigl\{h \ov {\lambda(r^*, p)}\bigr\} \cdot\exp\{c \ov{\mu(r^*, p)}\},
\quad \text{where $r^*(z):=\ov{r(\ov {z^{-1}})}^{-1}$,   }
\end{equation}
determines a Hilbert space of holomorphic functions invariant with respect to 
the operators $\rho^{h,c}(\frQ)$ defined by {\rm(\ref{eq:sections})}.
The corresponding representation of Virasoro algebra is $L(h,c)$.
Moreover, for $q\in\Diff$ these operators are unitary.
 \end{corollary}

\section{Preliminaries. Explicit formulas for representations of $\Gamma$.}

\COUNTERS

{\bf\punct The Fock space.} See \cite{Ner-book}, Sect. V.3, or \cite{Ner-gauss}, Sect. 4.1, for details of definition and proofs.
The {\it boson Fock space} with $n$ degrees of freedom 
is a Hilbert space of holomorphic functions on $\C^n$ with inner
product 
$$
\la f, g\ra=\frac1{\pi^n}\int_{\C^n} f(z)\ov{g(z)} e^{-|z|^2}\,\prod_{j=1}^n d\Re z_j\, d\Im z_j.
$$
 For $a=(a_1,\dots,a_n)\in\C^n$ we define an element $\phi_a^{(n)}\in F_n$ by
 $$
\phi_a(z):=\exp\Bigl\{\sum_{j=1}^n z_j \ov a_j\Bigr\}.
$$
Then we have a {\it reproducing property}
$$
f(a)=\la f, \phi_a\ra_{F_n}.
$$

Consider an isometric embedding $J_n:F_n\to F_{n+1}$ determined
by
$$
Jf(z_1,\dots,z_n,z_{n+1})=f(z_1,\dots,z_n).
$$
Next, consider the union of the chain of Hilbert spaces
$$
\dots \longrightarrow F_n \longrightarrow F_{n+1} \longrightarrow\dots
\, .$$
Its completion is the {\it Fock  space $F_{\infty}$ with infinite number degrees of
freedom}.

For $a=(a_1,a_2,\dots)\in\ell_2$ we define an element $\phi_a\in F$
by 
$$
\phi_a(z):=\exp\Bigl\{\sum_{j=1}^n z_j \ov a_j\Bigr\}:= \lim_{n\to\infty}  \phi_{(a_1,\dots,a_n)}(z_1,\dots,z_n),
$$
the limit is a limit in the sense of the Hilbert space $F_\infty$.
Then for any element $h\in F_\infty$
we define a holomorphic function on $\ell_2$ by
$$
f(u)=\la f, \phi_u \ra_{F_\infty}.
$$
We identify $F_\infty$ with this space of holomorphic functions.

\sm

Since any infinite-dimensional Hilbert space $H$ is isomorphic $\ell_2$, we can define a Fock space
$F(H)$ as a space of functions on $H$ determined by the reproducing kernel
$$
R(h_1,h_2)=\exp\, \la h_1,h_2\ra_H.
$$
On Hilbert spaces determined by reproducing kernels, see, e.g., \cite{Ner-gauss}, Sect. 7.1.

\sm

{\bf\punct Symbols of operators.} For details, see \cite{Ner-book}, Sect. V.3, VI.1,
\cite{Ner-gauss}, Sect.4.1, Sect.7.1. Let $A$ be a bounded linear operator in $F_\infty$.
Following F.A.Berezin, define its {\it kernel}
by
$$
K(z,u)=\la A \phi_u,\phi_z\ra_{F_\infty}.
$$
Then for any function $f\in F_\infty$, we have
$$
f(z)=\la f, k_z\ra_{F_\infty}.
$$
where $k_z(u):=\ov{K(z,u)}$.

Let $A$, $B$ be bounded operators, let $L$, $K$ be their kernels,
let $k_z(u):=\ov{K(z,u)}$, $l_z(u):=\ov{L(z,u)}$ then
the kernel $M$ of $BA$ is
$$
M(z,u)=\la k_u , l_z \ra_{F_\infty}.
$$

\sm

{\bf\punct Gaussian operators.} See \cite{Ner-book}, Sect. V.4, Sect. VI. 2-4,
 \cite{Ner-gauss}, Sect.5-6.
Consider a block symmetric matrix%
\footnote{The symbol $L^t$ denotes the transposed matrix.}
$S=\begin{pmatrix}
    K&L\\L^t&M
   \end{pmatrix}
   $ of size $\infty+\infty$
   and two vectors $\lambda$, $\mu\in \ell_2$.
   Define a {\it Gaussian operator}
   $$
   B\begin{bmatrix}
    K&L&\bigr|&\lambda^t\\L^t&M&\bigr|&\mu^t
   \end{bmatrix}
   $$
in $F_\infty$ as an operator with the kernel
\begin{equation}
\exp\left\{\frac 12 \begin{pmatrix}
                     z & \ov u
                    \end{pmatrix}
   \begin{pmatrix}
    K&L\\L^t&M
   \end{pmatrix}   
   \begin{pmatrix}
                     z^t\\\ov u^t
                    \end{pmatrix}
                    + z \lambda^t+ \ov u\mu^t
                    \right\}.
                    \label{eq:gauss-expression}
\end{equation}
Here the vectors $z$, $u$, $\lambda$, $\mu\in\ell_2$ are regarded as vectors-rows,
$z^t$, $u^t$, $\lambda^t$, $\mu^t$ are vector-columns. The expression in the curly
brackets is a quadratic expression in $z$, $\ov u$.

If an operator $B[\dots]$ is bounded, then 

\sm

1. $\|S\|\le 1$, $\|K\|<1$, $\|M\|<1$;

\sm

2. $K$, $M$ are Hilbert-Schmidt operators.

\sm

If also $\|S\|<1$, then the operator $B[S|\sigma]$ is bounded.

\sm

More on conditions of boundedness, see \cite{Ner-book}, Sect.V.2-4.

\sm

Product of two Gaussian operators is given by the formula
{\footnotesize
\begin{multline}
    B\begin{bmatrix}
    K&L&\bigr|&\lambda^t\\L^t&M&\bigr|&\mu^t
   \end{bmatrix}\,
      B\begin{bmatrix}
    P&Q&\bigr|&\pi^t\\Q^t&R&\bigr|&\kappa^t
   \end{bmatrix}
   =
   \sigma(M, P; \mu, \pi)\cdot
   \\
   B\!\!\begin{bmatrix}
      K+LP(1-MP)^{-1}L^t& \!\!\!\!\!\! L(1-PM)^{-1} Q&\!\!\!\!\!\!\bigr|&\!\!\!\!\!\! \lambda^t+L(1-PM)^{-1} (\pi^t+P\mu^t)
      \\
      Q^t(1-MP)^{-1} L & \!\!\!\!\!\! R+Q^t(1-MP)^{-1} MQ &\!\!\!\!\!\!\bigr|&\!\! \!\!\!\!\kappa^t+Q^t(1-MP)^{-1} (M\pi^t+\mu^t)
     \end{bmatrix}
     ,
     \label{eq:gauss-product}
\end{multline}
}
where the constant $\sigma(\dots)$ is given by
\begin{equation}
 \sigma(M, P; \mu, \pi) =\det\bigl[(1-MP)\bigr]^{-1/2}
 \exp\left\{ \frac 12 \begin{pmatrix}
                       \pi&\mu
                      \end{pmatrix}
                      \begin{pmatrix}
                       -P&1\\ 1 &-M
                      \end{pmatrix}^{-1}
\begin{pmatrix}
                       \pi^t\\\mu^t
                      \end{pmatrix}
 \right\}
 .
 \label{eq:gauss-const}
\end{equation}

{\sc Remark.} The conditions $\|M\|<1$, $\|P\|<1$
imply $\|MP\|<1$. Hence we can evaluate  $(1-MP)^{-1/2}$
as 
$$
(1-MP)^{-1/2}=1+\frac 12 MP-\frac 38 (MP)^2+\dots
$$
Next, $M$, $P$ are Hilbert--Schmidt operators.
Therefore $MP$ is a trace class operator. It is easy to show that the series 
for $-1+(1-MP)^{-1/2}$ converges in the trace norm,  therefore the determinant
in (\ref{eq:gauss-const})
is well-defined.
\hfill $\boxtimes$

\sm

{\sc Remark.}
The big matrix in the right-hand side 
of (\ref{eq:gauss-product}) admits a transparent geometric interpretation, see \cite{Ner-book},
Theorem 5.4.3, Sect. VI.5. \hfill $\boxtimes$

\sm

{\bf \punct Gaussian vectors.} Define Gaussian vectors
$b[P|\pi^t]$, where $P$ is a Hilbert-Schmidt symmetric matrix
with $\|P\|<1$ and $\pi\in\ell_2$
by
$$
b[P|\pi](z):=\exp\Bigl\{\frac 12 zPz^t+ z\pi^t\Bigr\}.
$$
 Then $b[P|\pi^t]\in F_\infty$.
Also,
\begin{multline}
     B\begin{bmatrix}
    K&L&\bigr|&\lambda^t\\L^t&M&\bigr|&\mu^t
   \end{bmatrix}
   b[P|\pi^t]= 
\sigma(M, P; \mu, \pi)\times \\ \times  b[  K+LP(1-MP)^{-1}L^t|\lambda^t+L(1-PM)^{-1} (\pi^t+P\mu^t)],
\end{multline}
where $\sigma(\dots)$ is given by (\ref{eq:gauss-const}). Notice that
  $$
  K+LP(1-MP)^{-1}L^t\qquad \text{and} \qquad \lambda^t+L(1-PM)^{-1} (\pi^t+P\mu^t)
  $$ also are elements of the  right-hand side
  of (\ref{eq:gauss-product}). Inner  products of vectors $b[\dots]$ 
  are given by
  \begin{equation}
 \la  b[K|\mu^t], b[P,\pi^t]\ra_{F_\infty}=\sigma(M,\ov P; \mu,\ov \pi^t)
 .
  \end{equation}
  
  {\sc Remark.} We have $b[O|0]=1$. Obviously, for any Gaussian operator $B[S|\sigma]$,
  we have
  \begin{equation}
  \la B[S|\sigma]\cdot 1, \cdot 1\ra=1.
  \label{eq:1-1=1}
  \end{equation}
  
  {\bf \punct The Weil representation.} Now consider the space $H=\ell_2\oplus \ell_2$
  and consider the group of bounded operators
  \begin{equation}
   \begin{pmatrix}   
   \Phi&\Psi\\
   \ov\Psi&\ov\Phi
  \end{pmatrix}
  \label{eq:Phi-Psi}
,\end{equation}
which are symplectic in the following sense:
$$
\begin{pmatrix}
   \Phi&\Psi\\
   \ov\Psi&\ov\Phi
  \end{pmatrix}
  \begin{pmatrix}
   0&1\\-1&0
  \end{pmatrix}
  \begin{pmatrix}
   \Phi&\Psi\\
   \ov\Psi&\ov\Phi
  \end{pmatrix}^t
  =  \begin{pmatrix}
   0&1\\-1&0
  \end{pmatrix}.
$$

We denote by $\Sp(2\infty,\R)$ the group of such matrices satisfying an additional condition:

\sm

--- $\Psi$ is a Hilbert--Schmidt operator.

\sm 

Next, we consider the group $\ASp(2\infty,\R)$ of affine transformations of $\ell_2\oplus\ell_2$ generated
by the group $\Sp(2\infty,\R)$  and shifts by vectors of the form 
$(h,\ov h)$. According \cite{Ber}, \S4, Theorem 3, the group $\ASp(2\infty,\R)$ has a standard projective unitary representation
in $F_\infty$ by unitary Gaussian operators.
Elements of $\Sp(2\infty,\R)$ act by operators
\begin{equation}
\tau\cdot B\begin{bmatrix}
   \ov \Psi \Phi^{-1}& \Phi^{t-1}&\bigr|&0\\
   \Phi^{-1}&-\Phi^{-1}\Psi&\bigr|&0
  \end{bmatrix}
  .
  \label{eq:berezin}
\end{equation}
The shifts $(h,\ov h)$
act by the operators
$$
B\begin{pmatrix}
  0&1&\bigr|&h
  \\1&0&\bigr|&-h
 \end{pmatrix}
 .
$$
Let $g\in \ASp(2n,\R)$, let $B[S|\sigma]$ be the corresponding Gaussian operator.
By the formula (\ref{eq:gauss-product}), the matrix $S$ depends only on the linear part of $g$.

\sm

{\bf \punct The Hilbert space $V$.}
  First, we consider a space $V^{smooth}$ of smooth functions
    $
    f(z)=\sum c_k z^k
    $
  on the circle $S^1$ defined up to an additive constant.
  Define the inner product in $V^{smooth}$
   by 
  $$
  \la z^k, z^k\ra_{V^{smooth}}= |k|,\qquad   \la z^k, z^l\ra_{V^{smooth}}=0\quad\text{if $l\ne k$.}
  $$
 Denote by $V$ the  completion of $V^{smooth}$ with respect to this inner product.
 Notice that all elements of $V$ are $L^2$-functions on $S^1$.
 Denote by $V_+$ (resp. $V_-$) the subspace consisting of series $\sum_{k>0} c_k z^k$
 (resp. $\sum_{k<0} c_k z^k$). These subspaces consist of functions
 admitting holomorphic continuations to the disks $D_+^\circ$ and $D_-^\circ$
 respectively. 
 
 We define a skew-symmetric bilinear form on
 $V$ by the formula
 \begin{equation}
 \{f,g\}=\int_{|z|=1} f(z)\,dg(z).
 \label{eq:bracket}
 \end{equation}
 The subspaces $V_+$, $V_-$ are isotropic and dual one two another with respect to this form.
 The projection operators to $V_\pm$  are given  by
 the Cauchy integral
 \begin{equation}
 P_\pm f(z)=\pm \lim_{\epsilon\to 0_\mp}\frac 1{2\pi i} \int_{|u|=1+ \epsilon} \frac{f(u)}{u-z}
 .
 \label{eq:cauchy}
\end{equation}

For a diffeomorphism $\gamma\in\Diff$ we define the linear operator in $V$ by
\begin{equation}
T(\gamma)f(z)= f\bigl(\gamma^{-1}(z)\bigr)
\label{eq:T}
.
\end{equation}
Evidentely, our form is $\Diff$-invariant:
\begin{equation}
\bigl\{ T(\gamma)f_1, T(\gamma)f_2 \bigr\}=\bigl\{ f_1, f_2 \bigr\}
\quad\text{or}\quad
\bigl\{ T(\gamma)f_1,f_2 \bigr\}=\bigl\{ f_1, T(\gamma)^{-1}f_2 \bigr\}
.
\label{eq:invariance}
\end{equation}

 \sm

 {\bf \punct Construction of highest weight representations of $\Diff$.}
  Consider the Fock space $F(V_-)$  corresponding to the Hilbert space $V_-$.
 We wish to construct projective representations of $\Diff$ in $F(V_-)$ corresponding to the standard
 boson realization  of representations of $\Vir$. Let us regard the circle as the quotient $\R/2\pi\Z$
 with coordinate $\phi$.
 Fix $\alpha$, $\beta\in\R$. For $\gamma\in\Diff$ we consider the affine transformation
 of
 $V$ given by the formula
 \begin{equation}
 T_{\alpha,\beta}(\gamma^{-1})f(\phi)=f\bigl(\gamma(\phi)\bigr)+\alpha\bigl( \gamma(\phi)-\phi\bigr)+\beta \ln \gamma'(\phi)
 .
 \label{eq:two-parameters}
 \end{equation}
 
 {\sc Remark.}
 Recall that $\phi$,  $\gamma(\phi)$ are elements of $\R/2\pi\Z$. We take an arbitrary continuous
 $\R$-valued branch
 of $\gamma(\phi)$. Then the function $\gamma(\phi)-\phi$ is  well defined up to  an additive constant and it is an alement of $V$ and we can multiply it by a constant $\alpha$. \hfill $\boxtimes$
 
 \sm
 
 The formula (\ref{eq:two-parameters}) determines an embedding of the group
  $\Diff$ to the group $\ASp(2\infty,\R)$.
  Restricting the Weil representation to $\Diff$ we get a unitary projective
  representation of $\Diff$.
 On the level of the Lie algebra we obtain the representation determined in Subsection 
 \ref{ss:boson-realization}. See \cite{Ner-book}, Sect. VII.2
 
 \sm
 
Next, we wish to write the matrix (\ref{eq:Phi-Psi}).
 Applying (\ref{eq:cauchy}) we represent the block $\Phi$ in (\ref{eq:Phi-Psi}) as 
 $$
 \Phi f(z)= \lim_{\epsilon\to 0_-}\frac 1{2\pi i}\int_{|z|=1+\epsilon} \frac {f(\gamma^{-1}(z))\,dz}{z-u}
 .
 $$
 Formally, the right-hand side is well-defined for real-analytic $f$ and real analytic 
 $\gamma\in\Diff$.
 But it makes sense for any distribution $f$ and any smooth $\gamma$. For instance,
 we can decompose $f(\gamma^{-1}(z))$ into the Fourier series $\sum_{k\in \Z} c_k z^k$
 and take 
 $$
 \Phi f(z):= \sum_{k>0} c_k z^k.
 $$
 In the same way, we  write an expression for $\Psi$.
 
 In the next subsection, we will express the operators $\Phi^{-1}$, $\ov\Psi \Phi^{-1}$ in
the terms of welding.

 \sm
 
   {\bf \punct  Formulas for action of $\Gamma$.} For details and proofs, see \cite{Ner-sbornik},
   \cite{Ner-book}, Sect. VII.4-5.
 Let $\alpha$, $\beta\in\C$.  
  The semigroup $\Gamma$ acts in the space $F(V_-)$ by  Gaussian operators
  $$
  \cN_{\alpha,\beta} (\frR)=
  B\begin{bmatrix}
    K(r_+)&L(r_+,r_-) &|& -(\beta+i\alpha)\, \ell_1^t(r_+)+\beta\, m_1^t(r_+)\\
    L(r_+,r_-)^t& M(r_-)&|& -(\beta+i\alpha)\, \ell_2^t(r_-)+\beta\, m_2^t(r_-)
   \end{bmatrix}
   ,
  $$
where operators 
$$
K:V_-\to V_+,\quad L:V_-\to V_-,\quad L^t:V_+\to V_+,\quad
M:V_+\to V_-
$$
 and vectors 
 $$\text{$\ell^t_1$, $m_1^t\in V_+$},\qquad \text{$\ell^t_2$, $m_2^t\in V_-$}$$
will be defined  now.

Consider a function $f\in V_-$. Then the function $f\circ (r_+)^{-1}$ is defined on the contour
$r_+(S^1)$. Let us decompose%
\footnote{This decomposition is given by the Sokhotski--Plemelj integral formula.}
 it as $F_1+F_2$, where $F_1$ is holomorphic
in the domain $r_+(D_+^0)$, and $F_2$ is holomorphic
in $\ov\C\setminus r_+(D_+)$. This decomposition is determined up to
an additive constant, $F_1\mapsto F_1+c$, $F_2\mapsto F_2-c$.
We set
$$
K f:=F_1\circ r_+, \qquad L f=F_2\circ r_-.
$$
In a similar way, we take a function $f\in V_+$, decompose
$$
f\circ (r_-)^{-1}\Bigr|_{r_-(S^1)}=F_1+F_2
,$$
where $F_1$ is holomorphic in $r_-(D_-^\circ)$ and $F_2$ is holomorphic in
 $\ov\C\setminus r_-(D_-)$.
We set
$$
Mf:= F_1\circ r_-,\qquad L^t f= F_2\circ r_+
$$

Recall that the subspaces $V_+$ and $V_-$ are dual one to another with respect
to the pairing (\ref{eq:bracket}). It can be shown that the operators $L:V_-\to V_-$
and $L^t:V_+\to V_+$ are dual one to another. Also, the operators 
$M: V_+\to V_-$, $K:V_-\to V_+$ are dual to themselves.

Finally, we set 
\begin{align}
 \ell^t_1(r_+)=\ln \frac{r_+(z)}{z},\qquad 
  \ell^t_2(r_-)=\ln \frac{r_-(z)}{z},
  \\
  m_1^t(r_+)=\ln r'_+(z),\qquad  m_2^t(r_-)=\ln r'_-(z).
\end{align}

{\sc Remark.} The operator $K(r_+)$ is the {\it Grunsky matrix} of the univalent
function $r_+$, it is a fundamental object of theory of univalent functions,
 see \cite{Gru}, \cite{Gol}, \cite{Duren}.
\hfill $\boxtimes$

\sm

It remains to explain the meaning
of the expression (\ref{eq:gauss-expression}).
For%
\footnote{We write $f^t$, $\ell^t$, \dots instead of $f$, $\ell$, \dots to symbolize that in our formulas functions correspond to
	vector-columns (we apply an operator to a function from the left).}
 $f^t_\pm\in V_\pm$,
we set
$$
\begin{pmatrix} f_-&f_+\end{pmatrix}
\begin{pmatrix} K&L\\L^t&M\end{pmatrix}
\begin{pmatrix} (f_-)^t\\(f_+)^t\end{pmatrix}:=
\{f^t_-,K f^t_-\}+ 
2\{f^t_+,Lf^t_-\}
+
\{f^t_+,M f^t_+\}
,$$
 where the form $\{\cdot,\cdot\}$ is defined by (\ref{eq:bracket}), also
$$
\begin{pmatrix} f_-&f_+\end{pmatrix}\begin{pmatrix} \ell_1^t\\ \ell_2^t \end{pmatrix}
:=\{f^t_-,\ell^t_1\}+\{f^t_+,\ell_2\},
$$
we omit the similar expression with $m_1$, $m_2$.

Also we note that under this normalization we have 
\begin{equation}
\pi(\cN_{\alpha, \beta} (\frR) 1)=1,
\label{eq:1-1-1}
\end{equation}
where $\pi$ is the projection operator in $F(V_-)$ to
the line $\C\cdot 1$.

\section{The calculation}

\COUNTERS

{\bf \punct Step 1.} Concider
$$\frR=(\ov\C, r_+, r_-), \quad\frP=(\ov\C,p_+,p_-)\in\Gamma$$
such that 
$$
r_+(0)=0,\quad r_-(\infty)=\infty, \quad  p_+(0)=0, \quad  p_-(\infty)=\infty.
$$
Next, we choose $r^+$, $p^-$ 
such that $(\ov\C, r^+,r_-)$, $(\ov\C,p_+,p^-)\in\Diff$ 
as it was discussed in Subsection \ref{ss:lambda-mu}.
Take the corresponding $\gamma_r$, $\gamma_p$, 
and represent $\gamma_r \gamma_p$ as $(\ov\C, q_+,q_-)$.

According  formulas (\ref{eq:gauss-product}),  (\ref{eq:gauss-const}),
we get
$$
\cN_{\alpha,\beta} (\frR) \,\cN_{\alpha,\beta} (\frP)
= \kappa_{\alpha,\beta} (\frR,\frP)\, \cN_{\alpha,\beta} (\frR\circ \frP),
$$
where the
 canonical cocycle $\kappa_{\alpha,\beta}$ is given by
\begin{multline}
\kappa_{\alpha, \beta}(\frR,\frP)=
\det[(1-MK)]^{-1/2}\times
\\
\times \exp
\biggl\{\frac12 \begin{pmatrix}-(i\alpha+\beta) \ell_1+\beta m_1& -(i\alpha+\beta) \ell_2+\beta m_2 \end{pmatrix}
\begin{pmatrix}
 -K&1\\ 1&-M
\end{pmatrix}^{-1}
\times \\\times
\begin{pmatrix} -(i\alpha+\beta)\, \ell^t_1+\beta\, m^t_1
\\ -(i\alpha+\beta)\, \ell^t_2+\beta\, m^t_2
\end{pmatrix} \biggr\}.
\label{eq:1}
\end{multline}
where
\begin{align*}
 &K:=K(p_+), \qquad M:=M(r_-), \\
 &\ell_1:=\ell_1(p_+),\quad
 \ell_2:=\ell_2(r_-),\quad m_1:=m_1(p_+),\quad m_2:=m_2(r_-).
\end{align*}
In particular, we see that the cocycle
$$\kappa_{\alpha,\beta}(\frR,\frP)=\kappa_{\alpha,\beta}(r_-,p_+)$$
 does not depend on $p_-$, $q_+$ and is holomorphic in $\alpha$, $\beta$.
 
 \sm
 
 Keeping in the mind  further manipulations we represent (\ref{eq:1}) in the form
 \begin{multline}
 \kappa_{\alpha, \beta}(r_-,p_+)=
 \det[(1-MK)]^{-1/2}
 \times\\\times
 \exp\biggl\{
-\frac 12\alpha^2 \begin{pmatrix}\ell_1&\ell_2\end{pmatrix} 
\begin{pmatrix}
 -K&1\\ 1&-M
 \end{pmatrix}^{-1}
 \begin{pmatrix}\ell^t_1\\\ell^t_2\end{pmatrix} 
 -\\-i\alpha\beta
 \begin{pmatrix}\ell_1&\ell_2\end{pmatrix} 
 \begin{pmatrix}
 -K&1\\ 1&-M
 \end{pmatrix}^{-1}
  \begin{pmatrix}-\ell^t_1+m_1^t\\-\ell^t_2+m_2^t\end{pmatrix} 
  +\\+
 \frac12 \beta^2 \begin{pmatrix}-\ell_1+m_1&-\ell_2+m_2\end{pmatrix} 
   \begin{pmatrix}
   -K&1\\ 1&-M
   \end{pmatrix}^{-1}
   \begin{pmatrix}-\ell^t_1+m_1^t\\-\ell^t_2+m_2^t\end{pmatrix} 
 \biggr\}.
 \label{eq:11}
 \end{multline}
 
 \sm

{\bf\punct Step 2.}
We also use the notation of Subsection \ref{ss:cocycles}. In particular, we 
have another notation for the canonical cocycles $\kappa^{h,c}$. By definition,
$$
\kappa_{\alpha, \beta}(r_-,p_+)=\kappa^{\frac 12 (\alpha^2+\beta^2), 1+12\beta^2}(r_-,p_+).
$$
Setting $\alpha=\beta=0$, we get
$$
\kappa_{0,0}(r_-,p_+)=\kappa^{0,1}(r_-,p_+)=\det\bigl[(1-K(p_+)\,M(r_-))^{-1/2}\bigr].
$$
Obviously%
\footnote{Let $(h_1,c_1)$, $(h_2,c_2)$ be in the domain of unitarity (\ref{eq:unitarity-1}).
Let $v_1$, $v_2$ be the highest vectors in $L(h_1,c_1)$, $L(h_2,c_2)$. 
The cyclic span of $v_1\otimes v_2\in L(h_1,c_1)\otimes L(h_2,c_2)$
is the module $L(h_1+h_1,c_1+c_2)$. This implies the statement for  $(h_1,c_1)$, $(h_2,c_2)$ 
satisfying (\ref{eq:unitarity-1}). It remains to apply the holomorphy.%
\label{f:}},
\begin{equation}
 \kappa^{h_1,c_1}(\frR,\frP)\, \kappa^{h_2,c_2}(\frR,\frP) =  \kappa^{h_1+h_2,c_1+c_2}(\frR,\frP).
 \label{eq:power}
\end{equation}
Keeping in the mind the holomorphy, we get that $\kappa^{h,c}$ has the form
$$
\kappa^{h,c}(r_-,p_+)=\exp\bigl\{ h\lambda (r_-,p_+)+ c \mu (r_-,p_+)\bigr\},
$$
where $\lambda(r_-,p_+)$, $\mu(r_-,p_+)$ are some functions.
In another notation,
\begin{multline}
\kappa_{\alpha,\beta}(r_-,p_+)=\exp\Bigl\{ \frac 12 (\alpha^2+\beta^2) \lambda (r_-,p_+)+ (1+12\beta^2)
 \mu (r_-,p_+)\Bigr\}
=\\
=\det\bigl[\bigl(1-K(p_+))M(r_-)\bigr)^{-1/2}\bigr]\times\\\times \exp\Bigl\{ \frac 12 (\alpha^2+\beta^2) \lambda (r_-,p_+)+ 12\beta^2 \mu (r_-,p_+)\Bigr\}
.
\label{eq:2}
\end{multline}

{\bf\punct Step 3.}
The exponentials in (\ref{eq:1}) and (\ref{eq:2}) must coincide.  This implies
that the term with $\alpha\beta$ in (\ref{eq:11}) is absent.
This means that
\begin{equation}
 \begin{pmatrix}\ell_1&\ell_2\end{pmatrix} 
 \begin{pmatrix}
 	-K&1\\ 1&-M
 \end{pmatrix}^{-1}
 \begin{pmatrix}\ell^t_1\\ \ell^t_2\end{pmatrix} 
 =
  \begin{pmatrix}\ell_1&\ell_2\end{pmatrix} 
  \begin{pmatrix}
  -K&1\\ 1&-M
  \end{pmatrix}^{-1}
  \begin{pmatrix}m_1^t\\m_2^t\end{pmatrix} .
  \label{eq:13}
 \end{equation}
Next, we transorm the coeffitient at $\beta^2$ in (\ref{eq:11}) using (\ref{eq:13}),
\begin{multline*}
   \frac12  \begin{pmatrix}-\ell_1+m_1&\ell_2+m_2\end{pmatrix} 
   \begin{pmatrix}
   -K&1\\ 1&-M
   \end{pmatrix}^{-1}
   \begin{pmatrix}-\ell^t_1+m_1^t\\-\ell^t_2+m_2^t\end{pmatrix}
   =\\=
    \frac12  \begin{pmatrix}m_1&m_2\end{pmatrix} 
    \begin{pmatrix}
    -K&1\\ 1&-M
    \end{pmatrix}^{-1}
    \begin{pmatrix}m_1^t\\m_2^t\end{pmatrix} -
      \frac12  \begin{pmatrix}\ell_1&\ell_2\end{pmatrix} 
      \begin{pmatrix}
      -K&1\\ 1&-M
      \end{pmatrix}^{-1}
      \begin{pmatrix}\ell^t_1\\\ell^t_2\end{pmatrix}
\end{multline*}
Therefore, we can reduce (\ref{eq:11}) to the form
 \begin{multline}
 \kappa_{\alpha, \beta}(r_-,p_+)=
 \det[(1-MK)]^{-1/2}
 \times\\\times
 \exp\biggl\{
 -\frac 12(\alpha^2+\beta^2) \begin{pmatrix}\ell_1&\ell_2\end{pmatrix} 
 \begin{pmatrix}
 -K&1\\ 1&-M
 \end{pmatrix}^{-1}
 \begin{pmatrix}\ell^t_1\\\ell^t_2\end{pmatrix} 
 +\\+
 \frac 12 \beta^2 \begin{pmatrix}m_1&m_2\end{pmatrix} 
 \begin{pmatrix}
 -K&1\\ 1&-M
 \end{pmatrix}^{-1}
 \begin{pmatrix}m_1^t\\m_2^t\end{pmatrix}\biggr\}.
 \end{multline}
Comparing this with (\ref{eq:2}), we come to
\begin{align}
 \lambda (r_-,p_+)= - \begin{pmatrix}\ell_1&\ell_2 \end{pmatrix}
 \begin{pmatrix}
 -K&1\\ 1&-M
\end{pmatrix}^{-1}
 \begin{pmatrix}\ell_1^t\\\ell_2^t \end{pmatrix};
 \\
 \mu (r_-,p_+)= \frac 1{24}\begin{pmatrix}m_1&m_2 \end{pmatrix}
 \begin{pmatrix}
 -K&1\\ 1&-M
\end{pmatrix}^{-1}
 \begin{pmatrix}m_1^t\\m_2^t \end{pmatrix}.
\end{align}
Since
$$
 \begin{pmatrix}
 -K&1\\ 1&-M
\end{pmatrix}^{-1}=\begin{pmatrix}
M(1-KM)^{-1} & (1-MK)^{-1}\\
(1-KM)^{-1}& K(1-MK)^{-1}
\end{pmatrix},
$$
we have
\begin{equation}
 \lambda (r_-,p_+)=
 \ell_1\Bigl[ M(1-KM)^{-1} \ell_1^t+ (1-MK)^{-1} \ell_2^t \Bigr]+\qquad\qquad
 \label{eq:square}
 \end{equation}
 \begin{equation}
\qquad\qquad\qquad\qquad\qquad\qquad +\ell_2\Bigl\{ (1-KM)^{-1} \ell_1^t+ K(1-MK)^{-1} \ell_2^t \Bigr\}
 \label{eq:curly}
\end{equation}
Also we have a similar expression for $\mu(r_-,p_+)$.

\sm

{\bf\punct Step 4.}
Next, we evaluate  the following product by  formula (\ref{eq:gauss-product}) 
$$
\cN_{\alpha,\beta}(\gamma_r)\, \cN_{\alpha,\beta}(\gamma_p)=\mathrm{const}\cdot 
\cN_{\alpha,\beta}(\gamma_r\gamma_p).
$$
The formula involves 3 operators of the type $B[\dots|.]$. This implies
an identity involving 3 matrices $[\dots|.]$. We substitute 
$\beta=0$ and write  the expression for the last column
in the matrix $[\dots|.]$ corresponding to $\gamma_r\gamma_p$:
\begin{multline}
 \ell^t_1(q_+)= \ell^t_1(r^+)+ 
 \\+ L(r^+,r_-)
 \Bigl\{ \bigl(1-K(p_+) M(r_-)\bigr)^{-1} \bigl(\ell^t_1(p_+)+ K(p_+) \ell^t_2(r_-) \bigr)\Bigr\}
 ;
 \label{eq:curly-1}
\end{multline}
\begin{multline}
 \ell^t_2(q_-)= \ell^t_2(p^-)+\\+
 L^t(p_+,p^-)\Bigl[\bigl(1-M(r_-) K(p_+)\bigr)^{-1} 
 \bigl(\ell^t_2(r_-)+M(r_-) \ell^t_1(p_+)\bigr)\Bigr]
 \label{eq:square-1}
.
\end{multline}
Keeping in the mind identity
$$
M(1-KM)^{-1}=(1-MK)^{-1}M,
$$
we observe that the  expressions in curly brackets in (\ref{eq:curly-1}) and (\ref{eq:curly})
coincide. We express the term in the curly bracket from (\ref{eq:curly-1})
$$
\Bigl\{\dots \Bigr\}= 
 L(r^+,r_-)^{-1}\bigl(\ell^t_1(q_+)- \ell^t_1(r^+)\bigr)
$$
and substitute to (\ref{eq:curly}).
Also, the expressions in the square brackets in (\ref{eq:square-1}) and (\ref{eq:square})
coincide. We express the term in the square brackets from (\ref{eq:square-1})
and substitute it to (\ref{eq:square}). 
In this way, we get
\begin{multline}
\lambda(r_-,p_+) = \ell_2(r_-) L(r^+,r_-)^{-1} \bigl(\ell^t_1(q_+)-\ell^t_1(r^+) \bigr)
+\\+
\ell_1(p_+) L^t(p_+,p^-)^{-1}\left(\ell_2^t(q_-)-\ell^t_2(p^-) \right)
\label{eq:Lambda-transformed}
.
\end{multline}

{\bf\punct Step 5.} Now we wish to evaluate $L^{-1} (r^+,r_-)$.
For this aim consider the Gaussian operator, corresponding to the diffeomorphism $\gamma_r$,
$$
B\begin{bmatrix}
  K(r^+) & L(r^+,r_-) &\bigl|& *\\
  L^t(r^+,r_-) & M(r_-)&\bigl|& *
 \end{bmatrix}
$$
On the other hand  applying  formula
 (\ref{eq:berezin}) to the same $\gamma_r$ we get an expression of the form
 $$
 B\begin{bmatrix}
 	\ov \Psi \Phi^{-1}& \Phi^{t-1}&\bigr|&*\\
 	\Phi^{-1}&-\Phi^{-1}\Psi&\bigr|&*
 \end{bmatrix}.
 $$
 Therefore $L(r^+,r_-)=\Phi^{t-1}$, i.e. $L^{-1}=\Phi^t$.
Thus,
\begin{equation}
L(r^+,r_-)^{t-1}\ell_2^t=P_+ T(\gamma_r)\ell^t_2,
\label{eq:Lt1}
\end{equation}
where $P_+$ is the projection operators to $V_+$.
In the same way, we obtain
\begin{equation}
L(p_+,p^-)^{t-1}\ell^t_1=P_- T(\gamma_p) \ell^t_1.
\label{eq:Lt2}
\end{equation}

{\bf\punct Step 6.}
 Evidently, for 
$g\in V_-$, we have
\begin{equation}
\int_{|z|=1} g d(P_+ f)=\int_{|z|=1} g \,df= \bigl\{f,g\bigr\}.
\label{eq:last-1}
\end{equation}
In the same way, for  $f\in V_+$,
 $$
 \int_{|z|=1} g d(P_- f)=\int_{|z|=1} g \,df= \bigl\{f,g\bigr\}.
 $$
Now we can transform two summands in formula
(\ref{eq:Lambda-transformed}) for $\lambda(r_-,p_+)$. The first summand equals to
\begin{multline*}
\ell_2(r_-) L(r^+,r_-)^{-1} \bigl(\ell_1(q^t_+)-\ell^t_1(r^+) \bigr)=
\bigl(\ell_1(q_+)-\ell_1(r^+) \bigr)L(r^+,r_-)^{t-1} \ell^t_2(r_-)
=\\=\bigl(\ell_1(q_+)-\ell_1(r^+) \bigr)P_+ T(\gamma_r)  \ell^t_2(r_-)
=\\=
\Bigl\{\bigl(\ell_1^t(q_+)-\ell_1^t(r^+) \bigr)\,,\,P_+ T(\gamma_r) \ell^t_2(r_-)  \Bigr\}
=\\=
\Bigl\{\bigl(\ell_1^t(q_+)-\ell_1^t(r^+) \bigr)\,,\, T(\gamma_r) \ell^t_2(r_-)  \Bigr\}
=\\=
\Bigl\{T(\gamma_r)^{-1}\bigl(\ell_1^t(q_+)-\ell_1^t(r^+) \bigr)\,,\,  \ell^t_2(r_-)  \Bigr\}.
\end{multline*} 
We applied (\ref{eq:Lt1}), (\ref{eq:last-1}), and the invariance (\ref{eq:invariance}).

A similar calculation gives us the second summand of(\ref{eq:Lambda-transformed}),
\begin{multline*}
\ell_1(p_+) L^t(p_+,p_-)^{-1}\left(\ell^t_2(q_-)-\ell^t_2(p_-) \right)=\\=
\Bigl\{\ell^t_1(p_+) ,P_-T(\gamma_p)\left(\ell^t_2(q_-)-\ell_2^t(p_-) \right)\Bigr\}
=\\=
\Bigl\{\ell^t_1(p_+) ,T(\gamma_p)\left(\ell^t_2(q_-)-\ell_2^t(p_-) \right)\Bigr\}
=\\=
\Bigl\{T(\gamma_p)^{-1}\ell^t_1(p_+) ,\left(\ell^t_2(q_-)-\ell_2^t(p_-) \right)\Bigr\}.
\end{multline*}

This gives us formula (\ref{eq:lambda}) for $\lambda(r_-,p_+)$.

\sm

{\bf \punct Final remarks.}
The derivation of the formula (\ref{eq:mu}) $\mu(r_-,p_+)$ is similar, we simply change
$\ell_1$, $\ell_2$ to
$m_1$, $m_2$.

Our calculation of $\lambda(r_-,p_+)$, $\mu(r_-,p_+)$ is valid for $c\ge 1$, $h\ge (c-1)/24$. Due to the holomorphy we extend
the result to arbitrary $(h,c)\in\C^2$.

\sm

It remains to expain Corollaries 1.2--1.3 from the theorem.

\sm

{\bf\punct The action of $\Diff$ in the space of holomorphic functions.}
For an univalent function $s$, se denote $S^\circ:= \ov{s(\ov z)}$.
For any element $f$ of the Fock space $F(V_-)$ we assign a holomorphic functional $\Xi$
given by
$$
F(s):=\Bigl\la f, b\bigl[ K(s^\circ)|(-i \alpha+\beta) \ell(s^\circ)+ \beta m(s^\circ)\bigr] \Bigr\ra_{F(V_-)}
.$$
In this way, we send the Fock space to the space of holomrphic functionals on $\Xi$,
this gives us the desired realization (see \cite{Ner-sbornik}, Subs. 4.12).

\sm

{\bf \punct Reproducing kernels.} 
Let $c\ge 1$, $h\ge (c-1)/24$. The problem can be easily reduced to
 evaluation of inner products
\begin{equation}
\la \cN_{\alpha,\beta} (\gamma_1) \cdot 1\,,\, \cN_{\alpha,\beta} (\gamma_2)\cdot 1\ra_{F(V_-)}
\label{eq:77}
,\end{equation}
where $\gamma_1$, $\gamma_2\in\Diff$.
Since operators $\cN_{\alpha,\beta} (\gamma)$ are unitary up to a constant,
we have
$$\cN_{\alpha,\beta}(\gamma)^*=s\cdot\cN_{\alpha,\beta}(\gamma^{-1}).$$
for some $s\in\C$. On the other hand,
$$
\la \cN_{\alpha,\beta}(\gamma)\cdot 1, 1 \ra_{F(V_-)}= 
\la 1, \cN_{\alpha,\beta}(\gamma^{-1}) \cdot 1 \ra_{F(V_-)} =1
$$
and this implies $s=1$.
Therefore, (\ref{eq:77}) equals to
$$
\la \cN_{\alpha,\beta}(\gamma_2^{-1})\, \cN_{\alpha,\beta} (\gamma_1)  \cdot 1\,,\,  1\ra_{F(V_-)}
$$
and we reduce our problem to the evaluation of the canonical cocycle.
 
 \sm
 
Next, we have    apriory identity (see footnote \ref{f:})
for unitary representations $L(h_1,c_1)$, $L(h_2,c_2)$:
$$
K^{h_1,c_1}(r,p)\, K^{h_2,c_2}(r,p)= K^{h_1+h_2,c_1+c_2}(r,p).
$$
Taking arbitrary $L(h_1,c_1)$, $c_2>1$, and sufficiently large $h_2$
we can obtain $K^{h_1,c_1}(r,p)$ as
$$
K^{h_1,c_1}(r,p) = \frac{K^{h_1+h_2,c_1+c_2}(r,p)} {K^{h_2,c_2}(r,p)}.
$$


\noindent
\tt Math.Dept., University of Vienna,
 \\
 Oskar-Morgenstern-Platz 1, 1090 Wien;
 \\
\& Institute for Theoretical and Experimental Physics (Moscow);
\\
\& MechMath.Dept., Moscow State University;
\\
\& Institute for Information Transmission (Moscow).
\\
e-mail: neretin(at) mccme.ru
\\
URL:www.mat.univie.ac.at/$\sim$neretin

\end{document}